\documentclass[10pt, twocolumn]{article}

\usepackage[left=2cm,right=2cm,top=2cm,bottom=2cm]{geometry}
\usepackage{lineno,hyperref}
\usepackage{graphicx}
\usepackage[font=footnotesize,labelfont=bf]{caption}
\usepackage{tabularx}
\usepackage{subfig}
\usepackage{epsfig} 
\usepackage{amsmath} 
\usepackage{amssymb}  
\usepackage{euscript}
\usepackage{color}
\usepackage{algorithmic}
\usepackage{algorithm}
\usepackage{bm}
\usepackage{enumerate}
\usepackage{bbm}

\def\R{{\mathbb{R}}}
\def\C{{\mathbb{C}}}
\def\Z{{\mathbb{Z}}}
\def\Exp{{\mathbb{E}}}

\def\V{{\mathcal V}}
\def\E{\mathcal{E}}

\def\G{\mathcal{G}}

\def\M{\mathcal{M}}

\def\N{\mathcal{N}}

\def\av{\mathbf{a}}
\def\bv{\mathbf{b}}

\def\vv{\mathbf{v}}
\def\yv{\mathbf{y}}
\def\sv{\mathbf{s}}
\def\iv{\mathbf{i}}
\def\uv{\mathbf{u}}
\def\xv{\mathbf{x}}
\def\wv{\mathbf{w}}
\def\pv{\mathbf{p}}
\def\qv{\mathbf{q}}

\def\thetav{\boldsymbol{\theta}}

\def\alphav{\boldsymbol{\alpha}}
\def\betav{\boldsymbol{\beta}}

\def\xiv{\boldsymbol{\xi}}
\def\diag{\mathrm{diag}}
\def\1{\mathbbm{1}}

\def\figurescale{0.9}

\newcommand{\trace}{\mathrm{Tr}}
\newcommand{\armse}{\mathrm{ARMSE}}

\newcommand{\until}[1]{\{1,\dots, #1\}}
\newcommand{\fromto}[2]{\{#1,\dots, #2\}}

\newcommand{\proj}[1]{\langle #1\rangle}
\newcommand{\real}[1]{\mathrm{Re}\ #1}
\newcommand{\imag}[1]{\mathrm{Im}\ #1}

\newtheorem{theorem}{Theorem}

\newtheorem{proposition}[theorem]{Proposition}

\newtheorem{remark}[theorem]{Remark}

\title{Smart Grid State Estimation with PMUs Time Synchronization Errors} 

\author{Marco Todescato\footnote{M.~Todescato is with the Bosch Center for Artificial Intelligence, Renningen, Germany. e-mail: mrc.todescato@gmail.com. The work was completed during his Postdoctoral fellowship at DEI, University of Padova, Italy.}, \ 
	Ruggero Carli,$^\dagger$ \ 
	Luca Schenato\footnote{R.~Carli and L.~Schenato are with the Department
of Information Engineering, University of Padova, via Gradenigo 6/b, 35131, Padova, Italy.
 [carlirug,schenato]@dei.unipd.it}, \ 
	Grazia Barchi\footnote{G.~Barchi is with the Institute for Renewable Energy, Eurac Research, viale Druso 1, 39100, Bolzano, Italy.
e-mail: Grazia.Barchi@eurac.edu}}
        	 
\begin{document}

\maketitle

\begin{abstract}
We consider the problem of PMU-based state estimation combining information coming from ubiquitous power demand time series and only a limited number of PMUs. Conversely to recent literature in which synchrophasor devices are often assumed perfectly synchronized with the Coordinated Universal Time (UTC), we explicitly consider the presence of time-synchronization errors in the measurements due to different non-ideal causes such as imperfect satellite localization and internal clock inaccuracy. We propose a recursive Kalman-based algorithm which allows for the explicit off-line computation of the expected performance and for the real-time compensation of possible frequency mismatches among different PMUs. Based on the IEEE C37.118.1 standard on PMUs, we test the proposed solution and compare it with alternative approaches on both synthetic data from the IEEE 123 node standard distribution feeder and real-field data from a small medium voltage distribution feeder located inside the EPFL campus in Lausanne.
\end{abstract}

\section{Introduction} 
Traditionally addressed at the transmission level of the power system where it is formulated as a nonlinear weighted least squares problem \cite{fs:1970}, by far, State Estimation (SE) is one fundamental task to properly operate the system \cite{aa-age:2004}. 
However, the \emph{smart-grid paradigm}, with increasing penetration of renewables and the transition from passive to active non-linear loads, e.g., electric vehicles and storage, while progressively filling the traditional separation between transmission and distribution networks, is calling for system-level solutions to address issues crossing the entire grid such as demand side management and demand-response \cite{pp-dd:2011} just as an example. This need, in turn, calls for flexible solutions which can be easily adapted and used across all the grid infrastructure, thus putting SE back in the spotlight both in transmission and distribution networks.\\
In this transitioning process, a central and potentially revolutioning role is played by Phasor Measurement Units (PMUs). \emph{Ideally} synchronized with the Coordinated Universal Time (UTC) \cite{Arbiter}, PMUs are able to provide the phasor of an electrical waveform, i.e. voltage or current, thus enabling the possibility of directly measuring rather than estimating the state variables \cite{agp-jst-kjk:1986}. 
While initially conceived for transmission systems, because of their accuracy and fast reporting rate, PMUs are acquiring significant interest also at the distribution level where can be used to support different applications -- protection and stability assessment \cite{df-cdb:2012}, power quality evaluation \cite{ac-nl-cm:2009}, management of fast time-varying loads \cite{ic-etal:2009} -- among which SE \cite{mp-etal:2015,lz-2017,XIA:2015} is one of the most relevant.\\
To avoid large phase errors and achieve high measurement accuracy, PMUs require almost prefect synchronization with the reference time (UTC). Hence they are usually equipped with expensive GPS units which, in turn, compromise their large-scale use. Consequently, given the limited number of PMUs that can be deployed in the system, counterbalanced by their high performance, most solutions to the SE problem propose to combine perfectly synchronized phasorial measurements, coming from a small number of PMUs, with measurements gathered from conventional smart meters \cite{ls:2014}. Conversely, different source of uncertainty are present in PMUs with time synchronization representing only one of them \cite{gb:2015}. Indeed, GPS provides 1pps (pulse-per-second) synchronization signals, with a theoretical accuracy of $1\mu$s which affects synchronization offsets, while internal PMUs clocks present frequency deviations which can produce large time skews.\\ 
In view of the this, the SE problem in the presence of synchronization error has been studied in the literature. In \cite{sb-rc-mt:2014} a first investigation on the effect of sync error is provided together with a static distributed estimator suitable for distribution-only grids. In \cite{py-zt-aw-an:2013}, thanks to a small angles approximation, the authors formulate a measurement model which is bilinear w.r.t. grid state and sync error parameters. Then, two parallel Kalman filters are used to approximately and simultaneously solve for state and sync error estimation. The very same model is exploited also in \cite{jd-sm-ycw-hvp:2014}. Finally, in \cite{jab-ur-jc:2015} only PMUs offset error due to GPS is considered and a precision time protocol to support the synchronization is proposed.

In this paper, in the spirit of \cite{ls:2014} and building upon our preliminary work \cite{Todescato:2017},
we address the problem of SE based on PMU measurements and load pseudo-measurements
where we explicitly consider the presence of time synchronization error in the PMU measurements. In particular, we are interested in the simultaneous estimation of the grid state as well as of the synchronization error parameters (i.e., offset and skew). We show how such error easily leads to poor estimation performance if the measurement model does not properly account for it, even with high PMUs penetration.

The contribution of the paper is twofold: i) a Kalman-based algorithm for the simultaneous estimation of the state and of the synchronization error parameters in \emph{real-time}. Based on the linear model proposed in \cite{sb-fd:2015} we are able to approximate the power flow manifold around any feasible working point and to explicitly compute \emph{off-line} accurate expected performance without resorting to time-consuming Monte-Carlo simulations. Also, as a byproduct of the model choice, our methodology seamlessly applies to transmission as well as to distribution grids. ii) Resorting mainly to the IEEE C37.118.1 \cite{ieee-pmu-standard} standard on PMUs, we quantify the impact of synchronization errors in state estimation and show that, if not compensated for, they can shadow the possible benefits of using PMUs altogether. To complement our contributions, we compare the proposed solution against \cite{ls:2014} on synthetic data from the IEEE 123 node distribution feeder as well as against the recently proposed estimator \cite{mp-etal:2015, lz-2017} on real-field data collected from the smartgrid located inside the EPFL campus \cite{epfl_grid}.


\section{Network Model}\label{sec:network_model}
Here we present the linearized power network model used. Building upon the recent \cite{sb-fd:2015} to which we refer the interested reader for all the mathematical details, validation and assessment of the model,
we model an AC power network under synchronous sinusoidal steady-state condition as a graph $\G(\V,\E)$. The nodes set $\V=\until{n}$ denotes the electric buses, while the edges set $\E$ denote the set of electric branches between connected buses. In synchronous steady-state regime, for each bus $h\in\V$, we define:
\begin{itemize}
\item $u_h=v_h e^{j\theta_h}\in\C$ complex voltage at the bus terminal where $v_h,\theta_h\in\R$ are the modulus and phase of the complex phasor, respectively;
\item $i_h\in\C$ complex current injected at the bus;
\item $s_h = p_h + jq_h$ complex (apparent) power absorbed by the bus where $p_h,q_h\in\R$ are the active (real) and the reactive (imaginary) power, respectively.  
\end{itemize}
Also, we define the \emph{nodal admittance matrix} $Y\in\C^{n\times n}$ element-wise as
$$
[Y]_{hk} =
\begin{cases}
\begin{array}{lc}
y_h^{sh} + \sum_{\ell\neq h} y_{h\ell}\,,& \mathrm{if}\ k=h\,;\\
-y_{hk}\,, & \mathrm{otherwise}\,;
\end{array}
\end{cases}
$$
where $y_{hk}$ is the admittance of the electric line $(h,k)$ connecting bus $h$ with bus $k$, while $y_{h}^{sh}$ is the shunt admittance (admittance to ground) at bus $h$.\\
By conveniently collecting all the nodal quantities into vectors $\uv = [u_1,\ldots,u_n]^T$, $\iv = [i_1,\ldots,i_n]^T$, $\sv = [s_1,\ldots,s_n]^T$, Kirchhoff's law and the nodal power balance read as
$$
\iv = Y\uv\, ,\qquad \sv = \diag(\uv)\overline{\iv}\, .
$$
where $\overline{(\cdot)}$ denotes the complex conjugate operator and $\diag(\cdot)$ denotes the diagonal matrix with $ii$-th diagonal element equal to the $i$-th element of its vector argument. Finally, by combining the two above equations, one gets
\begin{equation}\label{eq:node_balance}
\sv = \diag(\uv)\overline{Y\uv}\, ,
\end{equation}
which represent $n$ complex equations, i.e., the \emph{power flow equations}, that must be satisfied by any feasible power flow.\\
At this point, we recall the main result of \cite{sb-fd:2015} which will let us linearize the nonlinear power flow equations \eqref{eq:node_balance} around any feasible point in the power flow manifold.\\
We start by defining the power network state vector as 
$$
\xiv:=[\vv^T,\thetav^T,\pv^T,\qv^T]^T
$$ 
where $\vv,\thetav,\pv,\qv\in\R^n$ are obtained stacking together the corresponding nodal quantities. Then, by expressing the complex Eq.~\eqref{eq:node_balance} in rectangular coordinates, it is possible to rewrite them in implicit form as $F(\xiv)=0$, $F: \R^{4n}\mapsto\R^{2n}$, and, in turn, implicitly define (Lemma 1 of \cite{sb-fd:2015}) the \emph{power flow manifold}
\begin{equation}\label{eq:power_flow_manifold}
\M := \{\xiv\ |\ F(\xiv)= 0\}\, .
\end{equation}

\begin{proposition}[Proposition 1 of \cite{sb-fd:2015}]\label{prop:best_linear_approx}
Let $\xiv^*\in\M$, i.e., 
$$
\xiv^*=\{[(\vv^*)^T,(\thetav^*)^T,(\pv^*)^T,(\qv^*)^T]^T\quad |\quad F(\xiv^*)=0\, .\}
$$ 
Then, the linear manifold tangent to $\M$ in $\xv^*$ is given by
\begin{equation}\label{eq:linerized_pfes}
A_{\xiv^*}(\xiv-\xiv^*) = 0\, ,
\end{equation}
where
$$
A_{\xiv^*} = 
\Big[
\underset{A_{\uv^*}}{\underbrace{\left(\proj{\diag\ \overline{Y\uv^*}} + \proj{\diag\ \uv^*}N\proj{Y}\right)R(\uv^*)}} \quad -I
\Big]\, ,
$$
$\uv^*\!:=\!\vv^*e^{j\thetav^*}$, $I$ is the identity matrix of suitable size and
\begin{align*}
&N := \begin{bmatrix}I & 0 \\ 0 & -I \end{bmatrix}\, ,\quad
\proj{A} = \begin{bmatrix} 
				\real{A} & -\imag{A} \\ 
				\imag{A} &  \real{A} 
			\end{bmatrix}\, ,\\
&R(\uv) := 	\begin{bmatrix} 
			\diag(\cos \thetav)  &  -\diag(\vv)\diag(\sin \thetav) \\
			\diag(\sin \thetav)  &   \diag(\vv)\diag(\cos \thetav)
			\end{bmatrix}\, .
\end{align*}
\hfill$\square$
\end{proposition}

Proposition~\ref{prop:best_linear_approx} conveniently states how to reconstruct the best linear approximant, i.e., the plane tangent to $\M$ at $\xiv^*$, of the power manifold $\M$ at the feasible point $\xiv^*$. Interestingly, assuming $A_{\uv^*}$ invertible\footnote{Assuming $A_{\uv^*}$ invertible is not restrictive in real power grids given the presence of node shunt admittances.}, it is possible to express the voltage deviations $\delta\vv := \vv - \vv^*$ and $\delta\thetav := \thetav - \thetav^*$ (in polar coordinates) as linear functions of the power deviations $\delta\pv := \pv - \pv^*$ and $\delta\qv := \qv - \qv^*$ (in rectangular coordinates)
\begin{equation}\label{eq:linear_model}
\begin{bmatrix}
\delta\vv \\ \delta\thetav 
\end{bmatrix}
= 
A_{\uv^*}^{-1}
\begin{bmatrix}
\delta\pv \\ \delta\qv
\end{bmatrix}\, .
\end{equation}
Interestingly, the implicit formulation of Eq.~\eqref{eq:power_flow_manifold} defines all the voltages and power injections that are compatible with the physics without assuming any a priori model for the network's buses such as the typical PQ, PV or slack buses. Hence, as stated in \cite{sb-fd:2015}, one strength of the presented linear formulation is that it holds for any admissible working point $\xiv^*\in\M$ and that it can generalize many previously presented linear approximation such as the \emph{Linear Coupled power flow} model \cite{svd-ssg-ycc:2015}, the \emph{DC power flow} model \cite{bs-jj-oa:2009} and the \emph{rectangular DC power flow} model \cite{sb-sz:2016}. We refer to \cite{sb-fd:2015} for all the details.

\section{Measurement models}\label{sec:meas_model}
In the spirit of \cite{ls:2014}, we assume to have at disposal two types of information: i) ubiquitous historical data series of active and reactive power demands, used to provide the estimator with a rough prior knowledge of the state; and ii) sparse real-time high-accuracy phasorial measurements, used to refine the estimate.

\subsection{Power Demand Time-Series}\label{subsec:meas_power}
The first source of information are historical data series of active and reactive power demands collected at each bus from \emph{low-cost largely-available} and \emph{low accurate} smart meters. Namely, at node $h\in\V$ and time $t\in\Z_+$ we have
\begin{equation}\label{eq:power_meas}
\begin{bmatrix} \widetilde{p}_h(t) \\ \widetilde{q}_h(t) \end{bmatrix} =  
\begin{bmatrix} p_h \\ q_h \end{bmatrix} + 
\begin{bmatrix} w_h^p(t) \\ w_h^q(t) \end{bmatrix}\,,\
\begin{bmatrix} w_h^p(t) \\ w_h^q(t) \end{bmatrix} \sim
\N(0, \Sigma^w)
\end{equation}
where 
$$
\Sigma^w = 
\begin{bmatrix}
\sigma_p^2|p_h|^2 & \eta\sigma_p\sigma_q|p_h||q_h|\\ \eta\sigma_p\sigma_q|p_h||q_h| &  \sigma_q^2|q_h|^2 
\end{bmatrix}
$$
being $p_h,q_h$ the nominal values, $\sigma_p=\sigma_q\approx 30-50\%$ \cite{rs-rr:2014} and $\eta\in[0,1]$. Also, we assume $\Exp[w_k^p(t)w_h^p(t)] = \Exp[w_k^q(t)w_h^q(t)] = \Exp[w_k^p(t)w_h^q(t)] = 0$ \cite{ls:2014,rs-rr:2014}.

\subsection{PMU Measurements}\label{subsec:meas_pmu}
The second source of information are phasor measurements coming from \emph{high-cost} and \emph{highly accurate} \emph{Phasor Measurement Units} which, because of their cost, are usually deployed only at a limited number of electric buses. 
We recall that, to provide high accuracy phasorial values, PMUs are equipped with a GPS module exploited for synchronization purposes. Because of this, it is a far common assumption to consider PMUs as perfectly synchronized, thus neglecting the impact of the synchronization uncertainty on the resulting measurements. However, even in the presence of GPS modules, PMUs can still suffer of lack of synchronicity due to, e.g., temporary occlusion of satellites \cite{Castello:2018}. In addition to this, within successive synchronization instants with the GPS module, usually providing 1pps (pulse-per-second) synchronization signal, PMUs exploit an internal oscillator as reference clock which, in turn, might cause additional synchronization uncertainty. Hence, depending on the type of GPS module and oscillator, different synchronization accuracy, directly proportional to the cost of these devices, can be achieved.\\
Ultimately, the measurements at bus $h$ at time $t\in\Z_+$ are
\begin{equation}\label{eq:pmu_meas}
\begin{split}
&\widetilde{v}_h(t)\!=\!v_h(t)\!+\!w_h^v(t),\ \quad w_h^v(t)\sim\N(0,\sigma^2_{\mathrm{pmu},v} |v_h|^2),\\
&\widetilde{\theta}_h(t)\!=\!\theta_h(t)\!+\!w_h^\theta(t)\!+\!d_h(t),\ w_h^\theta(t)\sim\N(0,\sigma^2_{\mathrm{pmu},\theta}),
\end{split}
\end{equation}
where  we set\footnote{According to the synchrophasor reference standard IEEE C37.118:2014a \cite{ieee-pmu-standard}, PMUs must guarantee a TVE $<1\%$ in steady-state conditions. However, it is known that, especially at the distribution level, higher accuracy is required due to the lower power flows and angle phase differences \cite{Borghetti:11}.}
$\sigma_{\mathrm{pmu},v} = 0.1\%$, $\sigma_{\mathrm{pmu},\theta}=10^{-3}\mathrm{[rad]}$,
and assume uncorrelated measurement noise within the same node and across different nodes, i.e., $\Exp[w_k^v(t)w_k^\theta(t)] = \Exp[w_k^v(t)w_h^v(t)] = \Exp[w_k^\theta(t)w_h^\theta(t)] = 0$. Indeed, $w^v$ is mainly due to sampling jitter and synchronization error while $w^v$  to the instrumentation amplitude noise.
Differently from standard Gaussian additive models, the additional term $d_h(t)$ in \eqref{eq:pmu_meas} represents the error with respect to the true universal time. Since this component mainly affects the angle measurements \cite{gb:2015}, we restrict our analysis to synchronization errors (also referred to as de-sync) affecting voltage\footnote{We restrict the analysis to voltage only since the linear model \cite{sb-fd:2015} conveniently expresses them in polar coordinates. Conversely the currents are expressed in rectangular coordinate thus introducing additional approximation errors due to the projection of the measurement noise from polar to rectangular coordinates.} angles only. In particular, we consider a clock error model which, within successive synchronization instants $kT,(k+1)T,\ldots$, being $T$ the GPS synchronization period, has the form
\begin{equation}\label{eq:de-sync_model}
d_h(t) = \beta_h + \alpha_h\frac{T}{M-1}\ t\, ,\quad t\in\fromto{0}{M-1}\, ,
\end{equation}
where $M$ is the number of PMU measurements collected within two successive synchronization instants and
$\beta_h,\alpha_h\in\R$ are an offset term due to GPS error and the clock skew due to the fact that the internal clock of the PMU in general does not oscillate at the reference frequency, respectively. 

\smallskip

\begin{figure}
\centering
\includegraphics[width=\figurescale\columnwidth]{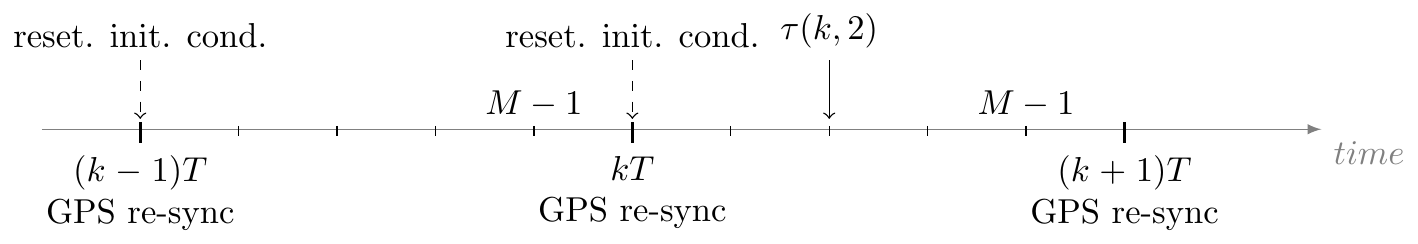}
\caption{Illustrative representation of the evolution of discrete time instants in a universal time frame. According to Algorithm~\ref{alg:estimate}, it is shown that every $T$[s] the filter reinitialize the initial condition.}
\label{fig:discrete_time_evol}
\end{figure}

\begin{remark}[Measurements time-scales]\label{rem:time_scales}
Observe that the measurements are characterized by three different time-scales. First, power demand information coming from historical time-series are usually available in the form of 1-day a-head predictions. 
Second, the PMU measurements live on a faster time-scale, $t,t+1,\ldots$, depending on the PMU reporting rate. 
Finally, as GPS provides a sync signal every $T$[s], the PMU internal clock re-synchronize with the universal reference at $kT,(k+1)T,\ldots$.
%
Then, it is convenient to denote discrete time instants in a universal time reference denoted as
$$
\tau(k,t) = kT + \frac{T}{M-1}t\, ,\quad k\in\Z_+\, ,\quad t\in\fromto{0}{M-1}\, ,
$$
referring to the $k$-th re-sync instant and to the $t$-th measurements within $[kT,(k+1)T)$, see Figure~\ref{fig:discrete_time_evol}. 
However, with a slight abuse of notation, in the following when talking about the evolution of a discrete quantity $x$, instead of writing $x(\tau(k,t))$ we use the simpler $x(k,t)$.\hfill$\square$
\end{remark}

\section{Estimation}\label{sec:estimation}
By taking advantage of the linear model \eqref{eq:linear_model}, the measurements and their statistical information, our final goal is to provide an accurate and real-time state estimate as output of a Kalman-based algorithm which, conversely to standard procedures, explicitly considers the effect of synchronization uncertainty in the measurements units. 
We cast our estimation procedure as a Bayesian inference process. More specifically, power demand predictions, characterized by a low accuracy, are used to provide a prior for our Bayesian model while highly accurate PMU measurements are used to improve the state estimate.

\subsection{State-Space Model}\label{subsec:power_model}
Our estimator consists of a Kalman filter \cite{rk:1960}, hence a suitable state-space model is needed.
Conventionally in the power system literature, since the network state are the voltages in rectangular coordinates, those are also chosen as filter state.
Conversely, as state vector for the filter, we pick the power demand deviations $\delta\pv,\delta\qv\in\R^n$ together with the synchronization error parameters\footnote{In general, the de-synchronization parameters equal the number $0\leq m \leq n$ of PMUs deployed in the grid.} $\alphav,\betav\in\R^m$. Our choice is not restrictive and naturally arises from \eqref{eq:linear_model} and \eqref{eq:pmu_meas}.
Now, since the incremental linear model \eqref{eq:linear_model} is defined w.r.t. a predefined operating point, assume $\pv^*$ and $\qv^*$ are nominal demands to which correspond $\uv^*=(\vv^*,\thetav^*)$. Then, state and output models at $\tau(k,t)$, $k\in\Z_+, t\in\fromto{0}{M-1}$, consists of
\begin{eqnarray}
&&\xv(k,t+1) = \xv(k,t) + \wv^{x}(k,t)\label{eq:ss_model}\\ 
&&\yv(k,t) = H\xv(k,t) + \wv^{y}(k,t)\label{eq:output_model}
\end{eqnarray}
\begin{eqnarray*}
&&\xv(k,t) = [\delta\pv(k,t)^T\ \delta\qv(k,t)^T\ \alphav(k,t)^T\ \betav(k,t)^T]^T,\\
&&\xv(k,0)\sim\N(0,\Sigma_0),\\
&&\Sigma_0 = 
\begin{bmatrix}
 \sigma_p^2\diag(|\pv^*|)^2 & \Sigma_0^{pq} & 0 & 0 \\ 
\Sigma_0^{qp} &  \sigma_q^2\diag(|\qv^*|)^2 & 0 & 0 \\
0 & 0 & \sigma_\alpha^2 I & 0 \\ 
0 & 0 & 0 & \sigma_\beta^2 I 
\end{bmatrix},\\
&&\Sigma_0^{pq} = \Sigma_0^{qp} = \sigma_p\sigma_q\diag(|\pv^*|)\diag(|\qv^*|),\\
%
%
&&\wv^{x}(k,t)\sim\N(0,W),\\
&&\yv(k,t) = 
\begin{bmatrix}
\delta\widetilde{\vv}(k,t) \\ \delta\widetilde{\thetav}(k,t)
\end{bmatrix}
:=
\begin{bmatrix}
\widetilde{\vv}(k,t) \\ \widetilde{\thetav}(k,t)
\end{bmatrix} - 
\begin{bmatrix}
\vv^* \\ \thetav^*
\end{bmatrix}\, ,\\
&&\widetilde{\vv}(k,t) = [\widetilde{v}_1(k,t)\ \hdots\ \widetilde{v}_m(k,t)]^T,\\ 
&&\widetilde{\thetav}(k,t) = [\widetilde{\theta}_1(k,t)\ \hdots\ \widetilde{\theta}_m(k,t)]^T,\\
&&H=
\left[ 
\begin{array}{c|c} 
A_{\uv^*}^{-1} & \begin{array}{cc} 0 & 0 \\ t\frac{T}{M-1}I & I \end{array}
\end{array}
\right],\\
&&\wv^{y}(k,t):=[\wv^v(k,t)^T\ \wv^\theta(k,t)^T]^T \sim\N(0,R),\\
&&R=\begin{bmatrix}
\sigma^2_{\mathrm{pmu},v} \diag(|\vv^*|)^2 & 0\\
0 & \sigma^2_{\mathrm{pmu},\theta} I
\end{bmatrix}.
\end{eqnarray*}
where $W$ can be used to embed information among the process noise covariance estimated, e.g., from data. In the following, given the small re-synchronization period, we usually consider $W=0$ and outline an interesting closed-form analysis. Since the above model is incremental with respect to the nominal value $\xiv^*$ and since the de-sync parameters can assume both positive and negative values, the state is reasonably initialized as a zero mean Gaussian random variable.

\begin{remark}[Exact linear output model]\label{rem:linear_desync_model}
It is worth observing that, due to the linear relation \eqref{eq:linear_model} between buses power $\pv,\qv$ expressed in rectangular coordinates and voltage $\vv,\thetav$ expressed in polar coordinates, the de-synchronization enters linearly in the output model \eqref{eq:output_model} without any further approximation. This is opposed to standard approaches in the literature where, to deal with linear models, the network state is expressed in rectangular coordinates, i.e., real and imaginary parts of the voltages. In this case, to resort to linear output models, the synchronization error must either be assumed or approximated as purely imaginary \cite{ls:2014,sb-rc-mt:2014}, under the additional assumption of small voltage angles differences.
\footnote{The assumption of small voltage angle differences usually holds for power distribution grids where the voltage values are clumped together in the proximity of the voltage value at the point of common coupling (PCC). However, the same does not hold in power transmission grids.\cite{borghetti2011synchronized}}.  
Also, even in the case when no sync error is considered, i.e., $d(t)=0$, observe that phasorial measurements are practically collected in polar coordinates. Hence, by expressing the output model with the same representation, we do not need any further manipulation of the data, i.e., projection from polar to rectangular coordinates, which, in turn, requires re-computation of the measurements correlation, inevitably introducing additional errors.\hfill$\square$
\end{remark}

\subsection{Synchronization-aware State Estimator}\label{subsec:kalman_filter}
Thanks to Eqs.~\eqref{eq:ss_model}--\eqref{eq:output_model}, we have at disposal a complete linear model to built a Kalman filter \cite{rk:1960} to simultaneously estimate network state and de-sync parameters. Algorithm~\ref{alg:estimate} describes what we refer to as \emph{Synchronization-aware State Estimator}, hereafter denoted with SASE. 
\begin{algorithm}[t]
\centering
\begin{algorithmic}
\REQUIRE $\Sigma_0$, $R$, $H$. Initialize $\Sigma(0) = \Sigma_0$.
\FOR[\texttt{//Offline}]{$t\in\fromto{0}{M-1}$}
\STATE Compute and store
\begin{eqnarray*}
&L(t+1) = (\Sigma(t)+W)H^T(H(\Sigma(t)+W)H^T + R)^{-1} \label{eq:klm_gain}\\
&\Sigma(t+1) = (I-L(t+1)H)(\Sigma(t)+W) \label{eq:klm_cov_update}
\end{eqnarray*}
\ENDFOR
\FOR[\texttt{//Online}]{$k\in\Z_+$}
\STATE Initialize $\widehat{\xv}(k,0)=0$
\FOR{$t\in\fromto{0}{M-1}$}
\STATE 
$
\widehat{\xv}(k,t+1) = \widehat{\xv}(k,t) + L(t+1)(\yv(k,t+1) - H\widehat{\xv}(k,t))
$
\ENDFOR 
\ENDFOR
\end{algorithmic}
\caption{SASE}
\label{alg:estimate}
\end{algorithm}
Observe that to run Algorithm~\ref{alg:estimate} values for $\pv^*$, $\qv^*$, $\vv^*$ and $\thetav^*$ are required from which $\yv$ and $H$ are derived. By leveraging the information coming from the available power demand time series, $\pv^*$ and $\qv^*$ are computed as one-day a-head predictions. Then, by means of a single full AC power flow computation it is possible to compute the corresponding values for $\vv^*$ and $\thetav^*$.
Note that, since $\Sigma(t)$ does not depend on the measurements, its evolution can be computed offline and stored for $t\in\fromto{0}{M-1}$ thus alleviating the computational burden. 
Finally, thanks to Eq.~\eqref{eq:linear_model}, the estimated voltages are equal to
$$
\begin{bmatrix}
\widehat{\vv} \\ \widehat{\thetav}  
\end{bmatrix}
= 
\begin{bmatrix}
\vv^* \\ \thetav^*  
\end{bmatrix}
+ 
A_{\uv^*}^{-1}
\begin{bmatrix}
\delta\widehat{\pv} \\ \delta\widehat{\qv}
\end{bmatrix}\, ,
$$
where, by partitioning $\Sigma$ as $\Sigma_0$, the covariance is given by
$$
\Sigma^{\uv} := 
\begin{bmatrix}
\Sigma^{v} & \Sigma^{v\theta}\\
\Sigma^{\theta v} & \Sigma^{\theta}
\end{bmatrix} 
= 
A_{\uv^*}^{-1}
\begin{bmatrix} 
\Sigma^{p} & \Sigma^{pq} \\
\Sigma^{qp} & \Sigma^{q} 
\end{bmatrix}
A_{\uv^*}^{-T}\, .
$$
As a side note, observe that both the model \eqref{eq:ss_model}--\eqref{eq:output_model} and Algorithm~\ref{alg:estimate} are outlined for $t\in\fromto{0}{M-1}$ for a given $k\in\Z_+$. As suggested by Figure~\ref{fig:discrete_time_evol}, since at $\tau(k,0)$ the PMUs re-synchronize with the GPS, the filter is re-initialized to reset $\alphav$ and $\betav$ and allow the computation of a new estimate. Similarly, newly available $\pv^*$, $\qv^*$ and new data can be used to recompute the model and $W$, respectively.

\section{Two-nodes case with no process noise ($W=0$)}\label{sec:toy_network}
Consider a network consisting of one load connected to one generator (the PCC, $v_{\rm pcc} = 1$, $\theta_{\rm pcc}=0$) through a purely inductive line with susceptance $b=-1$[p.u.], in the absence of shunt admittance. For the sake of the analysis, we assume the load is absorbing only active power $p$ while $q=0$. In this case the flat profile is a particular solution which can be chosen as linearization point. Thus, by leveraging the linear model \eqref{eq:linerized_pfes} one has $p = \theta$ and $v=1$ being $v$ and $\theta$ the voltage magnitude and the phase at the load, respectively. Notice that since we assumed $q=0$, the voltage magnitude is fixed and equal to $1$ thus only $p=\theta$ is of interest. Now, assume to collect, within two successive sync instants, $M$ phase measurements of the form \eqref{eq:pmu_meas} which, as in \eqref{eq:output_model}, can be expressed as
\begin{equation}\label{eq:toy_ex_meas_model}
\yv = 
\begin{bmatrix}
y_1\\ \vdots \\y_{M-1}
\end{bmatrix}
\!=\!
\begin{bmatrix}
1 & 1 & 0\\
1 & 1 & \frac{T}{M-1} \\
\vdots & \vdots & \\
1 & 1 & T
\end{bmatrix}\!\!
\begin{bmatrix}
\theta \\ \beta \\ \alpha
\end{bmatrix}\!+\! 
\begin{bmatrix}
w_1 \\ \vdots \\ w_{M-1}
\end{bmatrix}
\!=\!C\xv\!+\!\wv^y\,,
\end{equation}
with $\wv^y\sim\N(0,R)$, $R=\sigma_{\mathrm{pmu},\theta}^2I$, $\xv_0\sim\N(0,\Sigma_0)$, $\Sigma_0=\diag(\sigma^2_{\theta},\sigma^2_\beta,\sigma^2_\alpha)$. Furthermore for the sake of the analysis, let us assume absence of process noise, i.e., $w^x=0, W=0$ which, in the case of reasonably stable power demands within $T$[s], represents an acceptable first order approximation. Then, the posterior variance matrix in information form reads as
$$
\Sigma = \left(\Sigma_0^{-1} + C^TR^{-1}C\right)^{-1}
$$
and, thanks to \eqref{eq:toy_ex_meas_model}, after some tedious but straightforward algebraic manipulations, it is possible to compute $\Sigma=\Sigma(\sigma_{\mathrm{pmu},\theta},\sigma_\theta,\sigma_\beta,\sigma_\alpha,M,T)$ in closed form (reported in \cite{To-Ca-Sc-Ba:tech2017} for space reasons). Interestingly, it can be seen that, in the limit of the product $MT$, it holds that, 
\begin{equation*}\label{eq:toy_ex_limit_cov}
\underset{MT\to\infty}{\lim} \Sigma = 
\begin{bmatrix}
\frac{\sigma_\theta^2\sigma_\beta^2}{\sigma_\theta^2+\sigma_\beta^2} & -\frac{\sigma_\theta^2\sigma_\beta^2}{\sigma_\theta^2+\sigma_\beta^2} & 0\\
-\frac{\sigma_\theta^2\sigma_\beta^2}{\sigma_\theta^2+\sigma_\beta^2} & \frac{\sigma_\theta^2\sigma_\beta^2}{\sigma_\theta^2+\sigma_\beta^2} & 0\\
0 &  0 & 0
\end{bmatrix}\,.
\end{equation*}
and, in particular, as shown in \cite{To-Ca-Sc-Ba:tech2017}, that
\begin{equation}\label{eq:toy_ex_alphabeta_cov}
[\Sigma]_{22} = \sigma_{22}(\sigma_\theta\sigma_\beta)\,,\qquad [\Sigma]_{33} = \sigma_{33}\left(\frac{1}{MT^2}\right)\,.
\end{equation}
Hence, while for growing $M$ or $T$, $\sigma_{33}\to 0$ meaning that the uncertainty on the skew parameter goes to zero and, consequently, the parameter is perfectly estimated, residual uncertainty remains on both $\theta$ and $\beta$ for which $\sigma_{11}, \sigma_{22}\not\to 0$. As can be seen from the output matrix $C$, this is due to the fact that $\theta$ and $\beta$ are linearly dependent. Nonetheless, similarly to $\sigma_{22}$, even $\sigma_{11}$ and $\sigma_{12}$ are functions of the product $\sigma_\theta\sigma_\beta$. Thus $\sigma_{11},\ \sigma_{12},\ \sigma_{22}\to 0$ for $\sigma_\theta\sigma_\beta\to0$, meaning that if $\sigma_\beta=0$ then $\theta$ is perfectly estimated and viceversa. As highlighted later in the simulation section, this suggests that the different performance between our proposed SASE and what will be referred to as \emph{Ground Truth} (GT) is majorly due to this linear dependence.

\section{Simulations}\label{sec:simulations}
In this section we test the proposed SASE algorithm on two different data sets: i) synthetic data generated from the standard IEEE 123 nodes test bed \cite{Kersting2001}; ii) field-data collected from the smartgrid located inside the EPFL campus, Switzerland \cite{epfl_grid}. We use the Matlab Matpower package \cite{matpower} for power flow computations. Finally, if not differently specified, Table~\ref{tab:parameters} summarizes the parameters' values. Some observation are in order. First, regarding the PMU reporting rate, since it depends on the network frequency, here we consider a set of values for $M$. Second, given the relatively small sync period $T=1$[s], desync parameters $\beta, \alpha$ are assumed constant within the interval $[kT,(k+1)T)$. Third, for $\sigma_\alpha$ we assume PMUs equipped with quartz-crystal oscillator characterized by an accuracy $\approx 10\div 30$ ppm \cite{jv:1994}; also for $\sigma_\beta$  we assume a $50$Hz frequency signal with GPS $\approx 0.5\div 1\mu$s accurate \cite{Arbiter}. Finally, given the small re-synchronization period of 1[s], we assume no process noise, $W=0$.

\subsection{Synthetic data set - IEEE 123 nodes grid}\label{subsec:sim_ieee123}
We compare the SASE algorithm against: i) an online iterative version of the \emph{Bayesian Linear State Estimaion} algorithm presented in \cite{ls:2014} (denoted as BLSE) assuming no synchronization error in the measurements; ii) a \emph{Ground Truth} (denoted as GT) strategy assuming perfect knowledge and compensation of the de-synchronization error.
\begin{table}[t]
\centering
\normalsize
\begin{tabularx}{\columnwidth}{c|X|c}
Parameter & \centering Value [units] & Ref.\\
\hline
\hline 
& & \\
$T$ 				& 1[s] - 1pps gps resync signal	& \cite{dh-du-vg-dn-dk-mk:2001,py-zt-aw-an:2013}\\
$M$ 				& \{20,25,30,50,60\}[samples] 	& \cite{ieee-pmu-standard}\\
$\sigma_p,\sigma_q$	& $50\%$ 	& \cite{rs-rr:2014}\\
$\sigma_{\mathrm{pmu},v},\sigma_{\mathrm{pmu},\theta}$		& $0.1\%$, $10^{-3}$[rad]	& \cite{ieee-pmu-standard}\\
$\sigma_\alpha$	& $10^{-2}$[rad]		& \cite{jv:1994}\\
$\sigma_\beta$	& $2\times 10^{-4}$ [rad] 	& \cite{ieee-pmu-standard}\\
$\vv^*,\thetav^*,\pv^*,\qv^*$	& power flow nominal solution & 
\end{tabularx}
\caption{Parameters used in the simulations}
\label{tab:parameters}
\end{table}
The estimation performance is measured in terms of \emph{Average} of \emph{Root Mean Square Error} which, given any two $n$-dimensional vectors $\av(t),\bv(t)$, possibly functions of time, we numerically approximate over $N=500$ Monte Carlo runs as

\begin{equation}\label{eq:empARMSE}
\widehat{\armse}(\av(t), \bv(t), t) = \sqrt{\frac{1}{N} \sum_{i=1}^N \frac{1}{n} \sum_{h=1}^n |a_h(t) - b_h(t)|^2}
\end{equation}
Observe that, under the assumptions of linear model and correct measurement noise statistic, the matrix $\Sigma(t)$ can be used to compute the empirical ARMSE as
\begin{equation}\label{eq:thARMSE}
\armse(t) = \sqrt{\frac{1}{n}\trace(\Sigma(t))}\,.
\end{equation}
Also, Eq.~\eqref{eq:thARMSE} can be leveraged to perform model qualification checking, indeed, $\widehat{\armse}(t)\approx\armse(t)$ means the non-linear measurements statistic is effectively captured by the linear filter built on the approximation.\footnote{Conversely to SASE and GT, 
since BLSE ignores the presence of de-sync and, thus ignores the true measurements statistic, $\Sigma(t)$ cannot directly be used in \eqref{eq:thARMSE}. Hence, for BLSE, $\Sigma(t)$ used in \eqref{eq:thARMSE} is computed resorting to a modified Riccati equation comprising the measurements error statistic.}

\begin{figure}[t]
\centering
\subfloat[][Voltage\label{subfig:voltage_armse_fixed_pmu_meas}]{
\includegraphics[width=\figurescale\columnwidth]{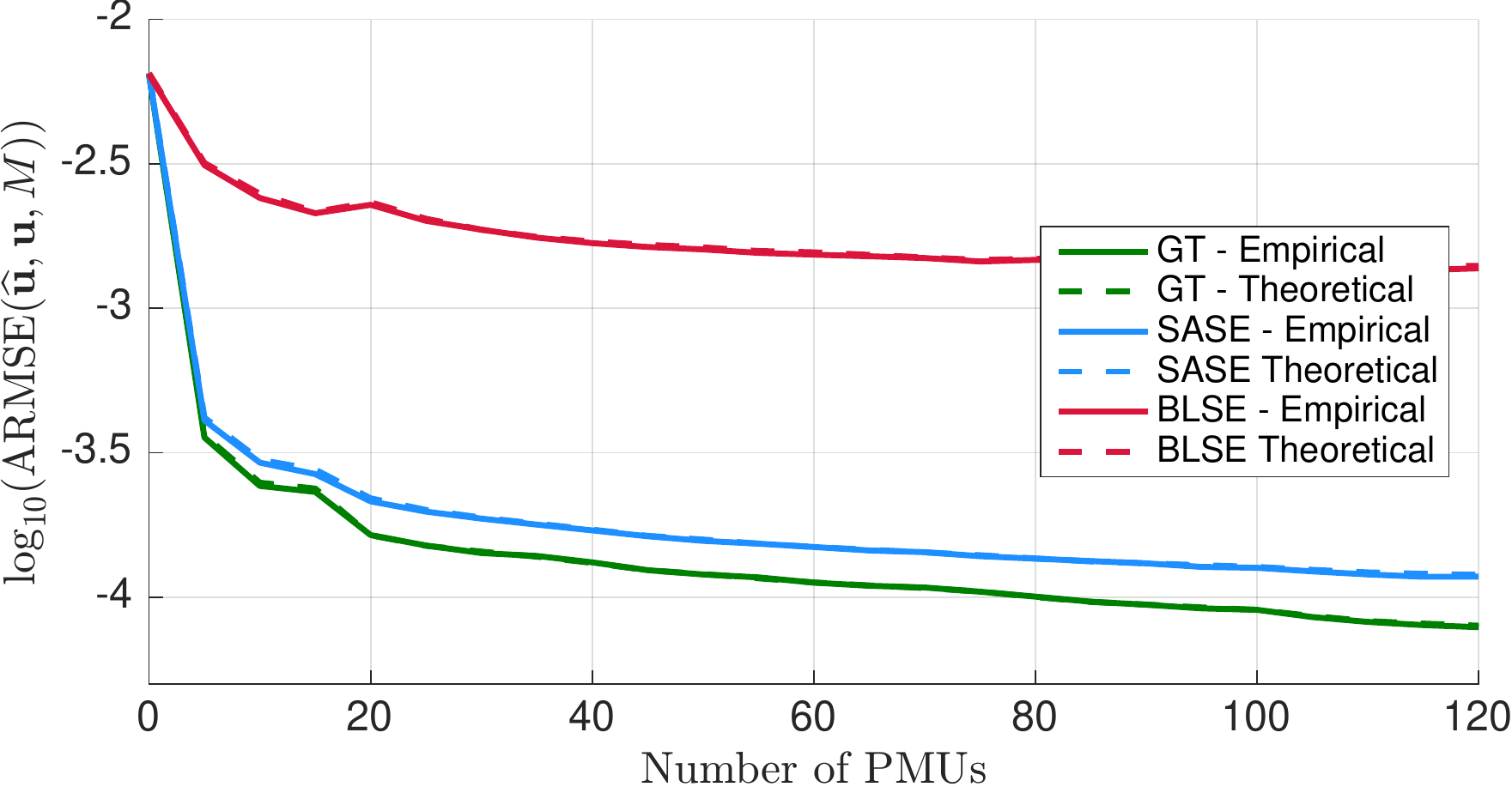}}\\
\subfloat[][Delay parameters\label{subfig:delay_armse_fixed_pmu_meas}]{
\includegraphics[width=\figurescale\columnwidth]{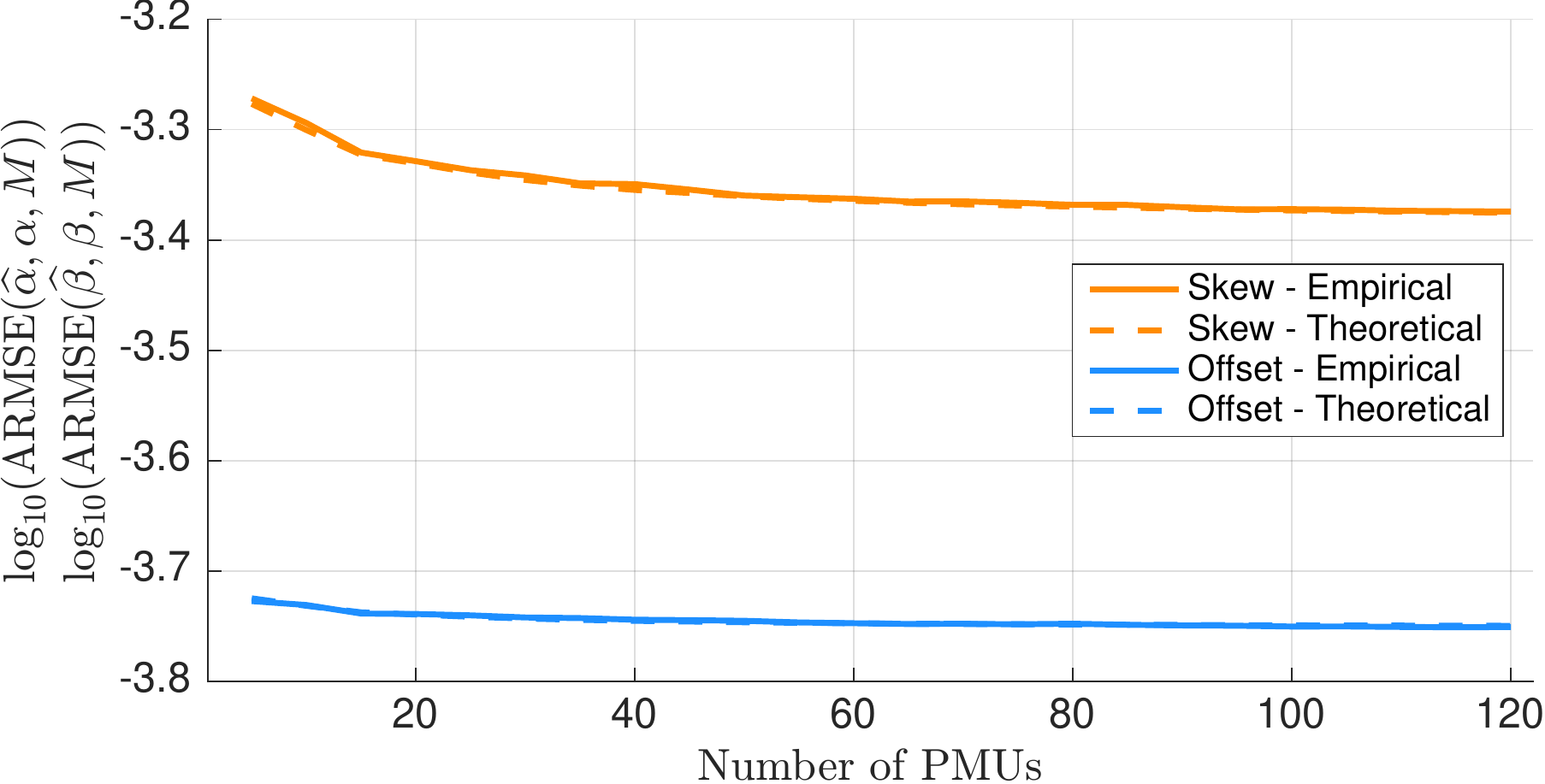}}
\caption{Empirical $\widehat{\armse}$ \eqref{eq:empARMSE} and theoretical $\armse$ \eqref{eq:thARMSE}, in logarithmic scale, as function of the number of deployed PMUs for $M=30$.}
\label{fig:armse_fixed_pmu_meas}
\end{figure}

%
Figure~\ref{fig:armse_fixed_pmu_meas} shows the performance as function of the number of PMUs deployed\footnote{For PMU deployment we exploitied the greedy approach presented in \cite{ls:2014}.} in the network when all the $M=30$ PMU measurements have been processed, i.e., right before a new synchronization instant occurs.
From Figure~\ref{subfig:voltage_armse_fixed_pmu_meas}, it is interesting to note that the proposed SASE approach behaves almost indistinguishable from the ground truth GT while the unmodeled synchronization error clearly deteriorates the BLSE, whose performance achieves, at best, $\approx 30\%$ improvement. Conversely, with only one PMU the proposed SASE performance improves of $\approx 60\%$.
Regarding the de-synchronization, Figure~\ref{subfig:delay_armse_fixed_pmu_meas} shows that the estimation performance does not improve for increasing number of PMUs deployed. This can be expected since the PMUs have been assumed uncorrelated. Finally, note that in both Figure~\ref{subfig:voltage_armse_fixed_pmu_meas} and \ref{subfig:delay_armse_fixed_pmu_meas}, theoretical and empirical curves almost perfectly coincide, validating the goodness of the prescribed linear approximation. This suggests that expensive Monte-Carlo simulations are not needed and many optimization problems such as optimal PMU placement or parameter sensitivity analysis can be performed very effectively also for large scale networks using the linearized model. It is worth stressing that this result holds for values of the parameters as in Table~\ref{tab:parameters}.\\
%
\begin{figure}[t!]
\centering
\subfloat[][Voltage\label{subfig:voltage_armse_fixed_num_pmu}]{
\includegraphics[width=\figurescale\columnwidth]{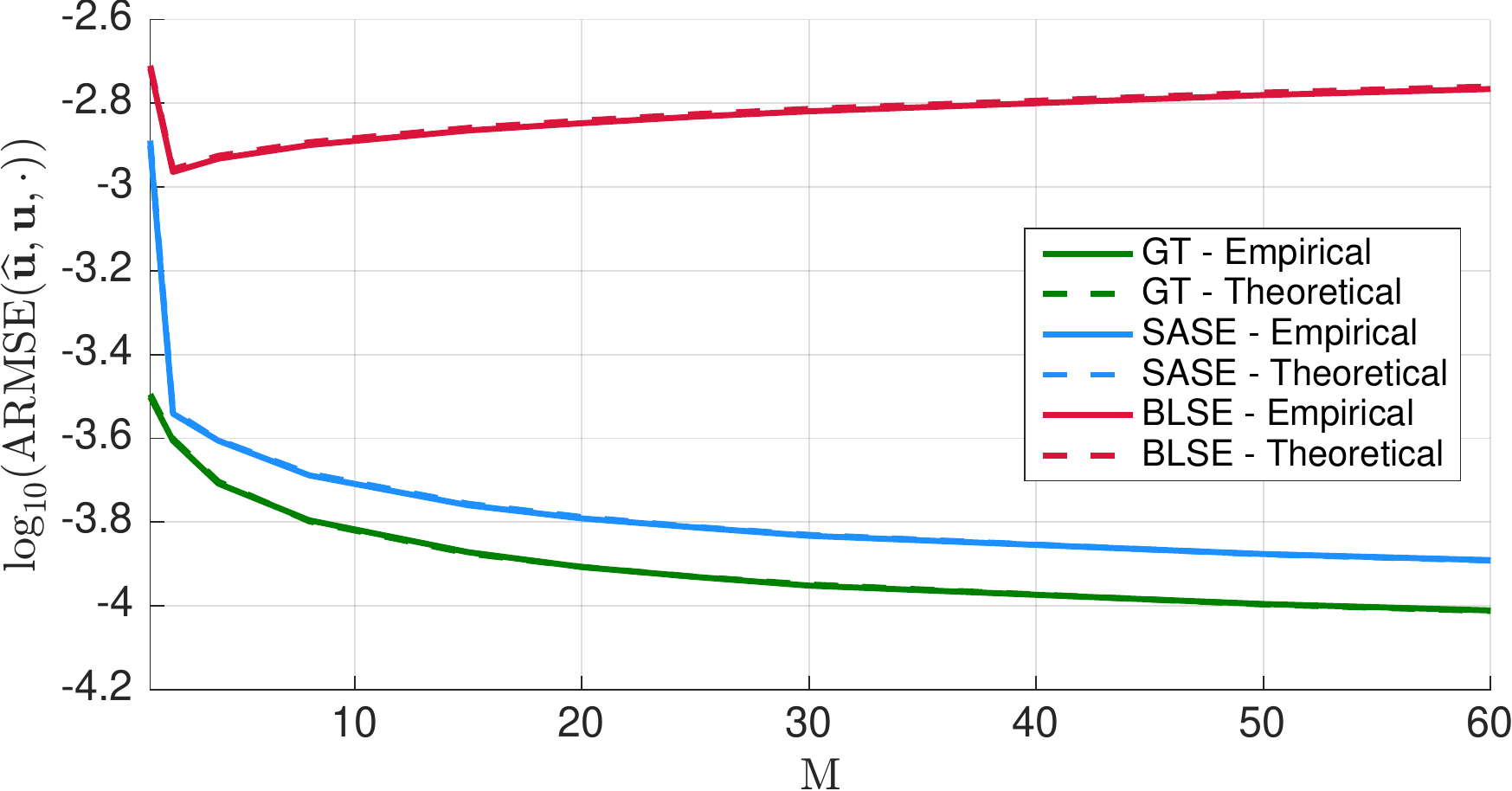}}\\
\subfloat[][Delay parameters. In dashed the theoretical values $\sigma_{22},\sigma_{33}$ of \eqref{eq:toy_ex_limit_cov} for offset and skew respectively, corresponding to the two-nodes case.\label{subfig:delay_armse_fixed_num_pmu}]{
\includegraphics[width=\figurescale\columnwidth]{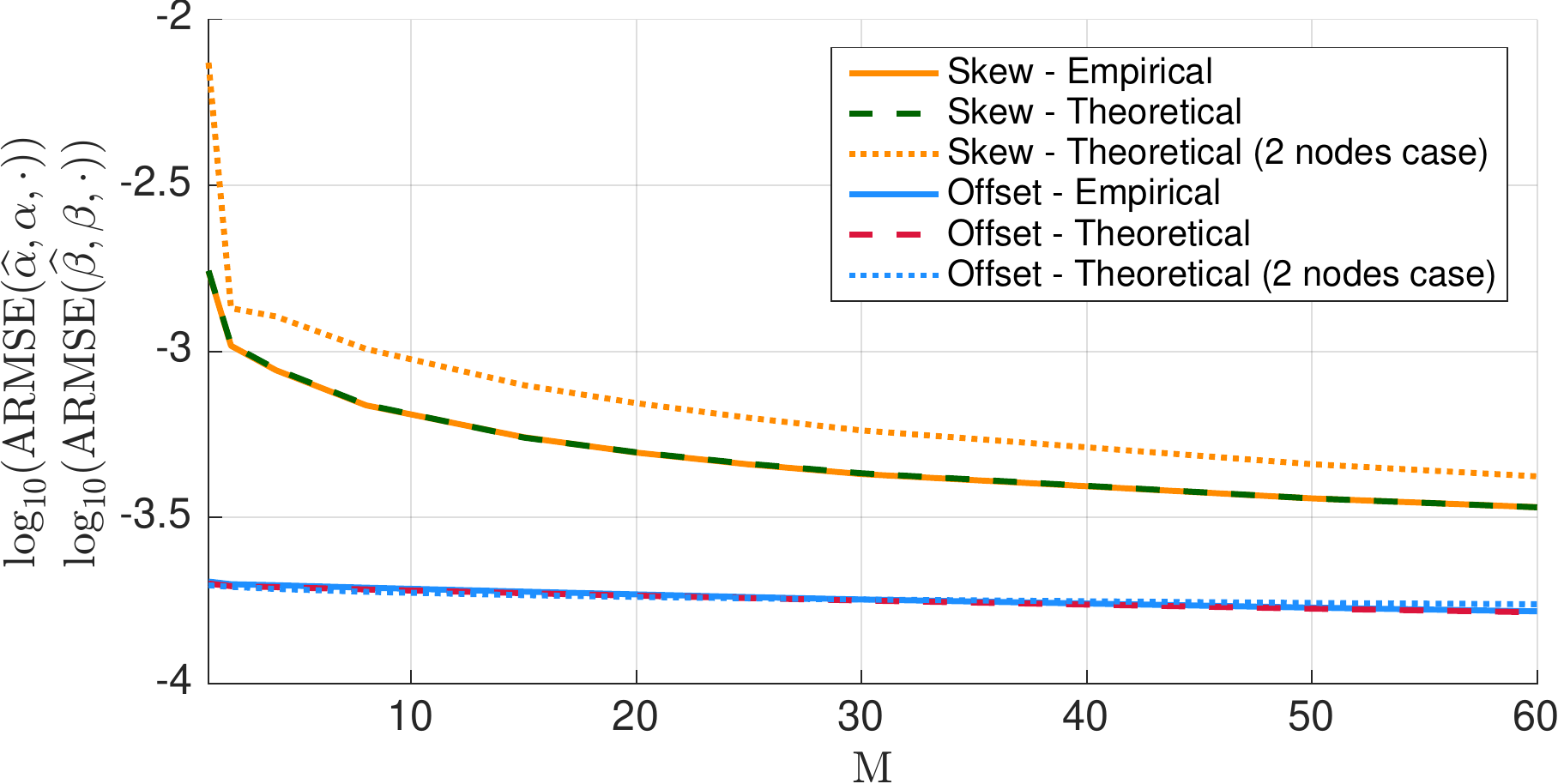}}
\caption{Empirical $\widehat{\armse}$ \eqref{eq:empARMSE} and theoretical $\armse$ \eqref{eq:thARMSE}, in logarithmic scale, as function of the number of collected PMU measurements $M$, for a fixed number of PMUs deployed in the network (in this case 8). $M\in[0,60]$ to consider the PMU reporting rates for both the $50$Hz and $60$Hz standards.}
\label{fig:armse_fixed_num_pmu}
\end{figure}
%
To emphasize the analysis of Section~\ref{sec:toy_network}, Figure~\ref{fig:armse_fixed_num_pmu} shows the performance for a fixed number of PMU as function of the number of collected PMU measurements $M$. Figure~\ref{subfig:voltage_armse_fixed_num_pmu} confirms the good behavior of the SASE compared with the GT. Conversely, as time passes, the BLSE does not improve its performance since it has no clue about the presence of the delay. Figure~\ref{subfig:delay_armse_fixed_num_pmu} supports this claim. Indeed, since the estimated skew improves for increasing $M$, the SASE is able to compensate for it. Finally, Figure~\ref{subfig:delay_armse_fixed_num_pmu} shows that the offset does not improve. As stressed in Section~\ref{sec:toy_network}, this is an intrinsic modeling problem due to the fact that offset and power demand happen to be linearly dependent. In addition, Figure~\ref{subfig:delay_armse_fixed_num_pmu} reports the values $\sigma_{22}$, $\sigma_{33}$ in \eqref{eq:toy_ex_alphabeta_cov} as a function of $M$ computed for the two-nodes network using the parameters value of Table~\ref{tab:parameters}. Observe how the theoretical values corresponding to the two-nodes case turn out to be extremely close to those obtained from the real network.  
This fact is interesting mainly for two reasons: i) it supports the claim that, in the limit for $M$ (or $T$), the proposed estimator perfectly reconstructs the skews while residual uncertainty remains in the offsets; ii) from the closed form expressions for $\sigma_{22}$ and $\sigma_{33}$, it is possible to retrieve, at least approximately, the value of the parameters needed to obtain a desired level of estimation accuracy.

\subsection{Real-world data set - EPFL smartgrid}\label{subsec:sim_epfl}
Here we test the SASE on data from the 20~kV 3-phase 6 nodes smartgrid installed in Lausanne within the framework of the NanoTera S$^3$-Grid project and located inside the EPFL campus \cite{mp-etal:2015,epfl_grid}. 
For a representation and an in-depth description of the network we refer to \cite{lz-2017}. We recall that the network is characterized by a line topology with nodes 1 to 5 monitored with PMUs (measuring current and voltage at 50Hz) and node 6 (the last along the line) is zero-injection and not monitored. Also, to estimate the measurements characteristics and noise variance values, we resort to the description reported in Section 4.2.1 of \cite{lz-2017} where the variances are computed from the data sheets of the PMU devices.\\
We compare the SASE against the first Kalman-based estimator proposed in \cite{mp-etal:2015, lz-2017} and carefully described in Sec 2.4 of \cite{lz-2017}.
A preliminary comparison between SASE and the algorithm in \cite{mp-etal:2015} highlights the following interesting facts:\\
i) The algorithm in \cite{mp-etal:2015} processes both currents and voltages in rectangular coordinates. This, even if the measurement covariances are correctly projected from polar to rectangular coordinates, inevitably introduces an approximation. SASE on the other hand considers only voltage measurements expressed in polar coordinates as naturally returned by the PMUs without introducing any additional source of error.\\
ii) The filter's state vector of \cite{mp-etal:2015} algorithm are real and imaginary parts of the voltages. Since voltage values are functions of all power demands, it is extremely hard to design meaningful priors, for both the mean and the covariance matrix, to properly initialize the filter. While this does not represent a critical issue in the presence of ubiquitous PMUs, in scenarios where the number of PMUs is small compared to the state dimension, this can translate in observability issues. Conversely, the SASE can be properly initialized leveraging, .e.g., one day-ahead forecast of power demands.\\
iii) Since the actual network frequency deviates from the nominal one, the phase angles are observed to rotate. To account for this rotation, at every iteration, the state model of the Kalman based algorithm in \cite{mp-etal:2015} is manually rotated of the quantity $\theta_k=2\pi\frac{f_k-f_0}{f_0}$, i.e., the phase angle difference due to the discrepancy between the nominal frequency $f_0$ and the actual frequency $f_k$ measured by the PMU itself. Hence, the algorithm resorts to additional information from the PMUs, namely, the measured network frequency or the frequency error (FE). Conversely, SASE automatically accounts for any linear trend existing in the phase angle measurements. \\
iv) Finally, in \cite{mp-etal:2015} no information regarding the posterior covariance characterizing the estimates are presented. Conversely, we show how the SASE nicely provides accurate estimates characterized by meaningful confidence intervals.

Before analyzing the simulation results we come back on the issue (see Section~\ref{subsec:power_model}) regarding the choice of the state model during synchronization instants $\tau(k,0)$, $k\in\Z$. Assume the state is $
\xv = [\delta\pv^T,\ \delta\qv^T,\ \alphav^T,\ \betav^T]^T
$ 
with evolution 
$$
\xv(k,t+1) = F_{\tau(k,t)} \xv(k,t) + \wv(k,t)\,. 
$$
To address the drift in the measurements, the state matrix is 
\begin{align*}
\begin{cases}
&F_{\tau(k,t)} = I\,,\quad
k\in\Z\,,\ t\in\fromto{1}{M-1}\,;\\
&F_{\tau(k,0)} = 
\begin{bmatrix}
I & 0 & 0 & 0\\
0 & I & 0 & 0\\
0 & 0 & I & 0\\
0 & 0 & I\cdot T & I
\end{bmatrix}\,,\quad
k\in\Z\,.
\end{cases}
\end{align*}
The last row-block of $F$ acts as an integrator for $\beta$ and is used to set their mean values at $\tau(k,0)$ to 
$$
\beta(k,0)=\beta(k-1,M-1)+\alpha(k-1,M-1)T.
$$

\begin{figure}[t]
\centering
\includegraphics[width=\figurescale\columnwidth]{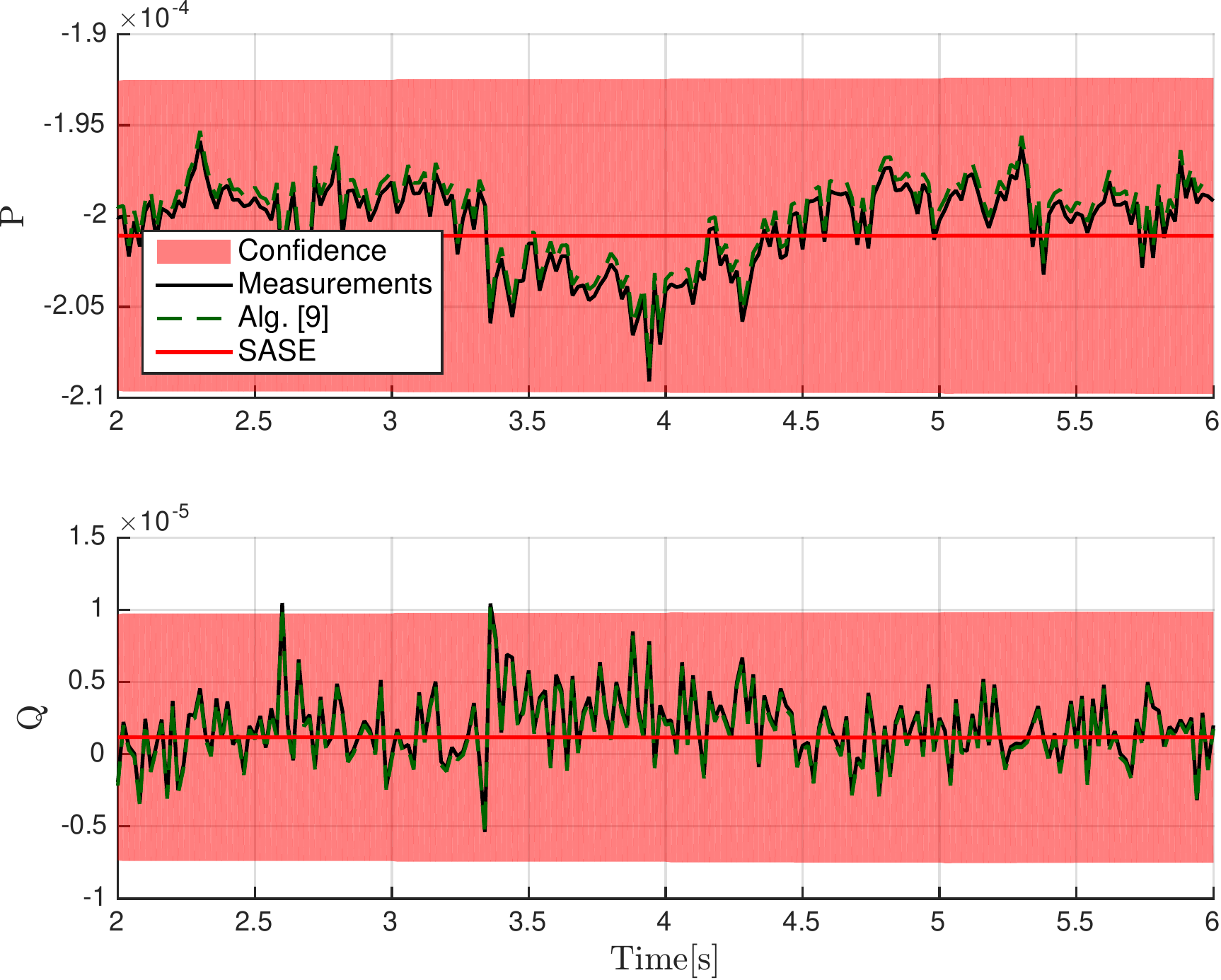}
\caption{Evolution of the estimates at node 3 (used for estimation) using three PMUs for estimation and two for validation.}
\label{fig:comparison_epfl_3pmu_estimate}
\end{figure}

\begin{figure}[t]
\centering
\includegraphics[width=\figurescale\columnwidth]{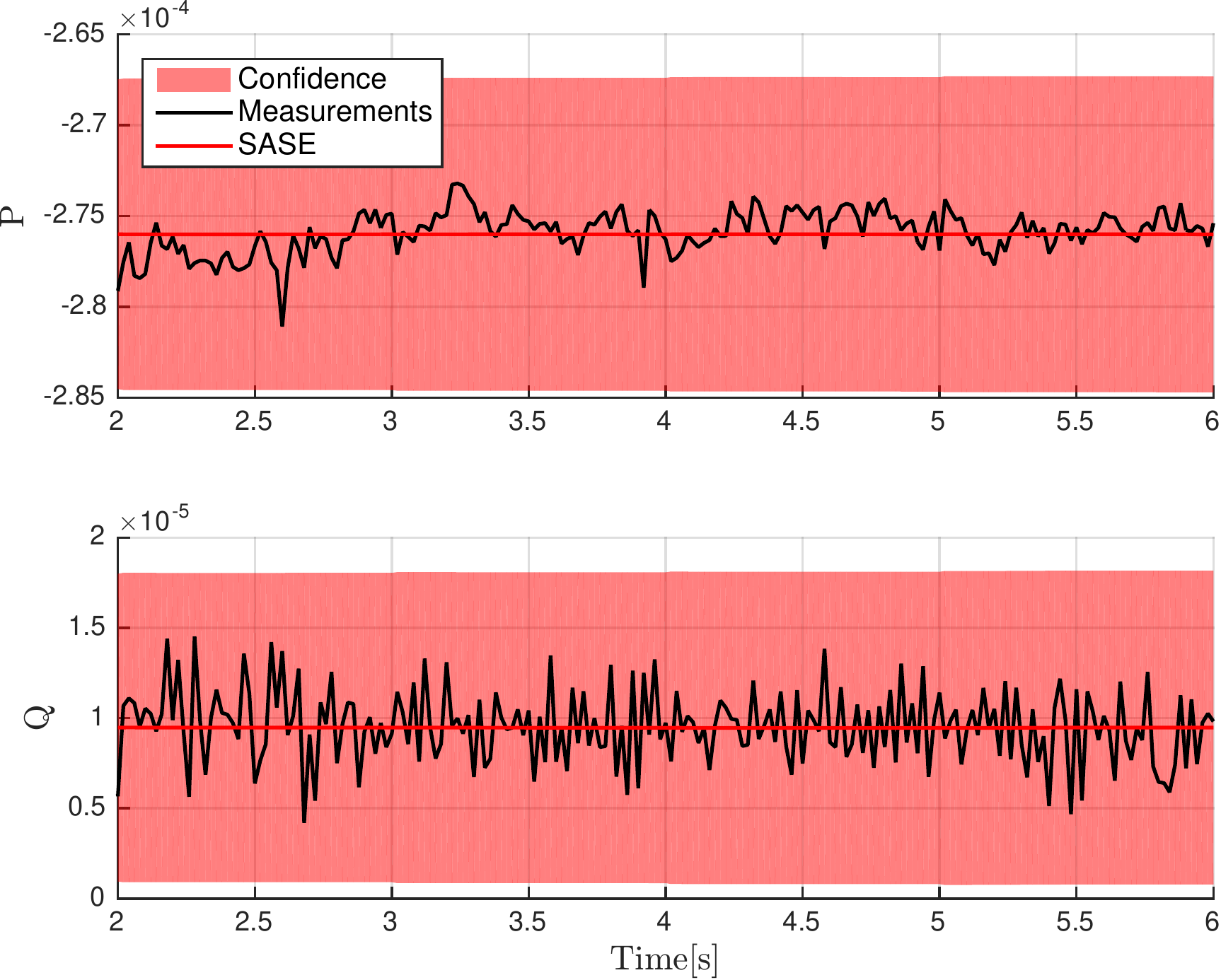}
\caption{Detail: Evolution of the estimates (excluding \cite{mp-etal:2015}) at node 5 (used for validation) using three PMUs for estimation and two for validation.}
\label{fig:comparison_epfl_3pmu_validation_noepfl}
\end{figure}

We now turn to the comparison between the two algorithms which are tested on a small subset of data consisting of a time window of 6 seconds collected on November 17$^{th}$, 2014, starting at 10:03:20~AM. More specifically, we are interested in comparing the SASE and the algorithm in \cite{mp-etal:2015} in terms of estimation and prediction. Hence, we assume to have at disposal measurements from nodes 1, 2 and 3 to perform the estimation while we use node 4 and 5 for validation, i.e., we do not use their measurements during estimation. The zero-injection node 6 is considered as virtual measurement for the algorithm in \cite{mp-etal:2015} and eliminated for the SASE.\\
Figures~\ref{fig:comparison_epfl_3pmu_estimate}--\ref{fig:comparison_epfl_3pmu_validation_noepfl} show the evolution of estimates and predictions at node 3 and 5, respectively, in the time interval $[2,6]$s. The first two seconds of simulation have been cut out in order to let the estimator in \cite{mp-etal:2015} to properly compute the covariance matrix $Q$ and converge to its steady state. From Figure~\ref{fig:comparison_epfl_3pmu_estimate} it is possible to see that the estimator in \cite{mp-etal:2015} nicely follows the measurements. The SASE, in harmony with its static state space model (since we assumed $W=0$), captures the average demand. However, the confidence interval returned by the estimator is in perfect accordance with the measurement values. Conversely, Figure~\ref{fig:comparison_epfl_3pmu_validation_noepfl} highlights the first difference between the two algorithms. Indeed, due to lack of (prior) information, in the prediction task, the algorithm in \cite{mp-etal:2015} does not provide any useful value (not even reported due to the high difference in the scaling factor). This analysis is supported by Table~\ref{tab:tve_3pmu} reporting the values of TVE of estimates and prediction w.r.t. the corresponding measurements.\\
%
\begin{table}[b]
\centering
\begin{tabular}{c|c|c|c}
Alg.\cite{mp-etal:2015} - est. & SASE - est. &  Alg.\cite{mp-etal:2015} - pred. & SASE - pred. \\ 
\hline 
\rule{0pt}{3ex}  $2.1\cdot 10^{-3}$ & $1.4\cdot 10^{-2}$ & $6\cdot 10^{5}$ & $8.5\cdot 10^{-3}$ 
\end{tabular}
\caption{Total vector error (TVE) for estimation and prediction of SASE and the algorithm \cite{mp-etal:2015} in the time interval $[2,6]$s, using 3 PMUs for estimation.}
\label{tab:tve_3pmu} 
\end{table}
%
In conclusion, even if the SASE algorithm was originally motivated to compensate de-synchronization with the GPS, it can be effectively used also to compensate linear frequency deviations at different nodes.

\section{Conclusions \& Future directions}\label{sec:conclusions}
In this paper we proposed a Kalman filtering procedure to address the problem of state estimation in power systems combining ubiquitous power time series and high-accurate sparse PMU measurements. Conversely to the standard assumption on the perfect synchronization of PMU devices, we considered the presence of PMUs de-synchronization and analyzed the problem of simultaneous estimation of network state and synchronization error parameters. Interestingly, by testing the proposed algorithm on both synthetic and real-field data, it is shown how the presence of synchronization errors can easily mislead the estimator if the measurement model does not properly account for it.

\bibliographystyle{IEEEtran}
\bibliography{IEEEabrv.bib,biblio.bib} %

\end{document}